\documentclass[11pt,notitlepage]{article}
\usepackage{graphicx}
\usepackage{amsfonts}
\usepackage{verbatim}
\usepackage{amsthm}
\usepackage{MnSymbol}
\usepackage{dsfont}
\usepackage{hyperref}
\usepackage{tabularx,array}
\usepackage{enumitem}
\usepackage{geometry}
\providecommand{\OO}[1]{\mathop{\mathrm{O}}\bigl(#1\bigr)}

\newtheorem{prop}{Proposition}
\newtheorem{theorem}[prop]{Theorem}
\newtheorem{lemma}[prop]{Lemma}

\newtheoremstyle{rem} 
{3pt}
{3pt}
{}
{}
{\bfseries}
{.}
{.5em}
{}
\theoremstyle{rem} 
\newtheorem{remark}{Remark}

\begin{document}

\noindent 

\begin{Large}
\begin{center}
{\bf Forgetting the starting distribution in finite interacting tempering}
 
 \vspace{1.2em}
Winfried Barta \\
\vspace{0.3em}
{\large \emph{George Washington University}} \\

\vspace{1.5em}
May 31, 2014

\end{center}
\end{Large}

\vspace{1em}

\begin{abstract}

Markov chain Monte Carlo (MCMC) methods are frequently used to approximately simulate high-dimensional,  multimodal probability distributions. In adaptive MCMC methods, the transition kernel is changed ``on the fly'' in the hope to speed up convergence. 
We study interacting tempering, an adaptive MCMC algorithm based on interacting Markov chains, that can be seen as a simplified version of the equi-energy sampler. Using a coupling argument, we show that under easy to verify assumptions on the target distribution (on a finite space), the interacting tempering process rapidly forgets its starting distribution. The result applies, among others, to exponential random graph models, the Ising and Potts models (in mean field or on a bounded degree graph), as well as (Edwards-Anderson) Ising spin glasses. As a cautionary note, we also exhibit an example of a target distribution for which the interacting tempering process rapidly forgets its starting distribution, but takes an exponential number of steps (in the dimension of the state space) to converge to its limiting distribution.
As a consequence, we argue that convergence diagnostics that are based on demonstrating that the process has forgotten its starting distribution might be of limited use for adaptive MCMC algorithms like interacting tempering.

\vspace{1em}

\end{abstract}

\vspace{1em}

\section{Introduction}
\label{sec:Intro}

Markov chain Monte Carlo (MCMC) algorithms are a widely used method to approximately sample from some complicated, often multi-modal probability distribution $\pi$ on a high-dimensional space $\mathcal{X}$. 
This is done by setting up a Markov chain $(X_t)$ that converges to $\pi$ as the number of steps $t$ goes to infinity. 
Many practical MCMC algorithms use local move Markov chains that can easily get ``stuck'' in one of the modes of the target distribution $\pi$. Tempering is a well-known strategy to try to overcome this problem. However, for some ``difficult'' distributions like the (ferromagnetic, mean-field) Potts model, even parallel and serial tempering algorithms are known to mix exponentially slowly in the dimension of the state space \cite{BhatnagarRandall2004, WoodardSchmidlerHuber2009b}.

The last decade has seen  considerable interest in \emph{adaptive} MCMC algorithms. Here, the transition kernel of the Markov chain depends on a parameter that may change over time, in a way that may depend on the entire history of the process so far. See \cite{AndrieuThoms2008, AtchadeFortMoulinesEtAl2011, Rosenthal2011} for recent overviews of these methods. 
\emph{Interacting tempering} \cite{FortMoulinesPriouret2011} is an adaptive MCMC algorithm based on several interacting Markov chains, each one targeting a tempered version of the distribution of interest $\pi$. It can be seen as a simplified version of the equi-energy sampler \cite{KouZhouWong2006}, which, in turn, attempts to improve on the convergence properties of the parallel tempering algorithm \cite{Geyer1991, Geyer2011}.
Since the interacting tempering process is generally not Markovian, standard Markov chain theory does not apply, and it takes considerable effort to establish 
ergodicity 
properties like convergence of marginal distributions, laws of large numbers, or central limit theorems. See \cite{AtchadeFortMoulinesEtAl2011,AtchadeWang2012,FortMoulinesPriouret2011, FortMoulinesPriouretEtAl2014, HuaKou2011, RobertsRosenthal2007} for important results in that direction. Quantitative, non-asymptotic rates of convergence are currently only poorly understood for adaptive MCMC algorithms in general, and interacting tempering in particular, but see \cite{PillaiSmith2013, SchmidlerWoodard2013} for first results in that direction.

In this work, we consider a  version of the interacting tempering algorithm for target distributions $\pi$ that satisfy two key requirements: First, the support of $\pi$ is \emph{simple} in the sense that it is easy to simulate the uniform distribution on it. Second, the distribution $\pi$ has exponentially bounded likelihood ratios, i.e. $\max_{x,y \in \mathcal{X}} \pi(x)/\pi(y) \leq \exp\{nD\}$ for some finite constant $D$ that does not depend on the dimension $n$ of the state space. Note that this implies that $\pi$ has bounded support. 
Examples of distributions that satisfy these assumptions are exponential random graph models, the Ising and Potts models (in mean field or on a bounded degree graph), as well as (Edwards-Anderson) Ising spin glasses, as discussed below. 
For ease of exposition, and since all examples  we have in mind are on finite spaces, we will assume that the state space $\mathcal{X}$ is finite throughout the paper.
Also note that when we speak of a probability distribution $\pi$ on a space $\mathcal{X}$, what we really have in mind 
is a family of distributions $(\pi^{(n)})_{n \in \mathbb{N}}$, where $\pi^{(n)}$ is a probability distribution on $\mathcal{X}_n$, and $n$ is the dimension of the state space $\mathcal{X}_n$. 
We are interested in the 
behavior of the algorithm as the dimension $n$ of the problem goes to infinity.

Our main result  is that the interacting tempering algorithm under these assumptions rapidly forgets its starting distribution (in order $n \log n$ steps). Importantly, since this process is not Markovian, our result says nothing about the (more interesting) question of how long it takes for the process to converge to its limiting distribution.
As a cautionary note, we exhibit an example of a distribution $\pi$ that satisfies these assumptions, 
but for which the interacting tempering algorithm takes an exponential (in $n$) number of steps to converge to its limiting distribution $\pi$.

In the absence of non-asymptotic, quantitative bounds on the convergence rates of MCMC algorithms, we often rely on \emph{convergence diagnostics}. A number of popular diagnostics work by demonstrating that the process in question has forgotten its starting distribution. 
However, for non Markovian processes like interacting tempering, forgetting the starting distribution is only necessary, but generally not sufficient, for convergence to the limiting distribution.
Our results illustrate that the time gap between 
forgetting the starting distribution and convergence to the limiting distribution
might often be huge, suggesting that diagnostics based on demonstrating forgetting of the starting distribution might be of limited use for adaptive MCMC algorithms like interacting tempering.

The rest of the paper is structured as follows: 
In section \ref{sec:Ass} we 
give a precise statement of our assumptions on the target distribution.
In section \ref{sec:examples}, we briefly discuss a number of well known models that satisfy these assumptions.
In section \ref{sec:Alg}  we give a precise definition of the interacting tempering algorithm in our setting. 
Section \ref{sec:Res} contains the statement and proof of our main result on rapid forgetting of the starting distribution for this algorithm. 
In section \ref{sec:NIH}, as a cautionary note, we exhibit an example of a target distribution for which our theorem implies rapid forgetting of the starting distribution, but where convergence to the limiting distribution takes at least an exponential number of steps. 
In section \ref{sec:Diag} we discuss  implications of our result for the use of convergence diagnostics.

\section{Assumptions}
\label{sec:Ass}
Throughout the paper we  assume that $\pi$ is a given probability distribution on a finite space $\mathcal{X}$. 
Without denoting this explicitly in the notation, we will always assume that there actually exists an entire family $(\pi^{(n)})_{n \in \mathbb{N}}$ of distributions, where $\pi^{(n)}$ lives on the state space $\mathcal{X}_n$ of dimension $n$.
Then, when we say that $\pi$ has a certain property, what we really mean by that is that each $\pi^{(n)}$ in the sequence satisfies the property in question. 
We are interested in what happens when the dimension $n$ of the problem goes to infinity.
Our central assumption on the distribution $\pi$ (really: on the sequence $(\pi^{(n)})_{n \in \mathbb{N}}$), is the following:
\begin{enumerate}[label=(\Alph*)]
	\item 	\begin{enumerate}[label=(\arabic*)]
			\item The support $\mathcal{X}$ of the distribution $\pi$ is  \emph{simple}, in the sense that the uniform distribution on $\mathcal{X}$  can be simulated in $\OO{n \log n}$ steps. 
			\\
			More specifically, we assume there exists a Markov transition kernel $P^{(0)}$ with unique stationary distribution $Uniform(\mathcal{X})$, that has the following property: For any $\epsilon>0$ there exists a constant $C(\epsilon)$, not depending on $n$, such that for any two starting states $v,w \in \mathcal{X}$ there exists a Markovian  coupling $(V_t,W_t)_{t \in \mathbb{N}_0}$ of two copies of the Markov chain $P^{(0)}$, started at $V_0=v$ respectively $W_0=w$, such that
			$$
			P_{v,w} \left\{ V_t \neq W_t \right\} \leq \frac{\epsilon}{n+1}, \;\;\;\;\; \mbox{ for all } t \geq C(\epsilon) n \log n.
			$$
			\item The distribution $\pi$ has exponentially bounded likelihood ratios, i.e. there exists  a constant $D$, not depending on $n$, such that
			$\max_{x,y \in \mathcal{X}} \pi(x)/\pi(y) \leq \exp\{nD\}$.
			\end{enumerate}
\end{enumerate}
Write $\pi_{max}$ and $\pi_{min}$ for the maximal and minimal values of $\pi(x)$ for $x \in \mathcal{X}$, respectively.
By defining $S(x) := n^{-1} \log \left( \pi(x) / \pi_{min} \right)$, we see that any probability distribution $\pi$ satisfying (A) can be written in the  form
%
%
%
%
\begin{equation}
\label{model}
\pi(x) := 
\frac{1}{Z(\beta)} \exp\{ n \, \beta \, S(x) \}, \;\;\;\;\;\; x \in \mathcal{X},
\end{equation}
where $\beta \geq 0$ is a constant, $S:\mathcal{X} \rightarrow [0,D]$ is a bounded non-negative function, and $Z=Z(\beta)$ is the normalizing constant.
Conversely, every probability distribution $\pi$ of the form (\ref{model}) satisfies assumption (A), if $\mathcal{X}$ is simple in the sense of (A)(1). 
This allows us to use representation (\ref{model}) in much of the rest of the paper.
On the other hand, the formulation in (A) is sometimes easier to check in applications, and we believe it better illustrates how mild the assumption is, and therefore how broadly applicable our results are.

The most serious restriction of our assumptions is finiteness of the state space $\mathcal{X}$. 
Our results could easily be extended to more general spaces, for example $\mathcal{X}=(0,1)^n$, as long as assumption (A) is satisfied. However, if the target distribution is specified via a density $\pi$ (with respect to some reference measure $\lambda$), it is clear that $(A)(2)$ implies that $\mathcal{X}$ has to be bounded, ruling out spaces like $\mathcal{X}=\mathbb{R}^n$. 
For ease of exposition, and since all applications we have in mind live on finite spaces, we restrict ourselves to the case of finite state spaces throughout this paper.

\section{Examples}
\label{sec:examples}
In this section we give some examples of probability distributions $\pi$ that satisfy our assumption (A) of
having simple support and exponentially  bounded likelihood ratios.

\vspace{0.5em}
\noindent
{\bf Ising models.}
The Ising model on a graph $G=(V,E)$ is the probability distribution
\begin{equation}
\label{eq:Ising}
\pi(\sigma) = \frac{1}{Z(\beta)} \exp \left\{ \beta \sum_{v \sim w} \sigma_v \sigma_w \right\}
\end{equation}
on $\mathcal{X} := \{-1,+1\}^V$, where the parameter $\beta \geq 0$ is called the inverse temperature, and $Z(\beta)$ is the normalizing constant (aka partition function). The sum in the exponent is over all edges 
$vw \in E$ 
of the graph $G$.
If the graph $G$ has bounded degree, and the number of vertices is $n$, then 
$$
2 \sum_{v \sim w} \sigma_v \sigma_w = \sum_{v \in V} \sigma_v \sum_{w: w \sim v} \sigma_w = \OO{n},
$$
so the model is of the form (\ref{model}).
We could also add an external magnetic field, resulting in
$$
\pi(\sigma) = \frac{1}{Z(\beta)} \exp \left\{ \beta \sum_{v \sim w} \sigma_v \sigma_w + h \sum_{v \in V} \sigma_v\right\}.
$$
Since this only adds an $\OO{n}$ term to the exponent, the model is still of the form (\ref{model}).
If $G=(V,E)$ is the complete graph on $n$ vertices, the temperature parameter is usually rescaled as $\beta := \alpha/n$ (to avoid trivial limits as $n$ goes to infinity), so that the model is of the form (\ref{model}) in this case as well.

\vspace{0.5em}
\noindent
{\bf Potts models.}
The Potts model with $q \geq 2$ colors, on the graph $G=(V,E)$ with $n$ vertices, is the probability distribution
\begin{equation}
\label{eq:Potts}
\pi(\sigma) = \frac{1}{Z(\beta)} \exp \left\{ \beta  \sum_{v \sim u} \mathds{1}_{\{ \sigma_v = \sigma_w \}} \right\}
\end{equation}
on $\mathcal{X} := [q]^V$, where $[q] := \{1,2,...,q\}$. Again, the parameter $\beta \geq 0$ is called the inverse temperature, $Z(\beta)$ is the normalizing constant, and the sum in the exponent is over all edges of the graph $G$.
If the graph $G$ has max  degree bounded by $d$, and we write $S(\sigma) := n^{-1} \sum_{u \sim v} \mathds{1}_{\{\sigma_u = \sigma_v\}}$, then we get
$$
\pi(\sigma) = \frac{1}{Z(\beta)} \exp \left\{ n \beta  S(\sigma) \right\},
$$
with $0 \leq S(\sigma) \leq d/2$. Thus, the distribution $\pi$ satisfies our assumption (\ref{model}) with $D := d/2$.
%
%

In mean field, when $G$ is the complete graph on $n$ vertices, we get the Curie Weiss Potts model
\begin{equation}
\label{eq:Potts}
\pi(\sigma) = \frac{1}{Z(\beta)} \exp \left\{ (\beta / n) \sum_{v,w \in V} \mathds{1}_{\{ \sigma_v = \sigma_w \}} \right\}
\end{equation}
on $\mathcal{X} := [q]^V$.
Since the sum in the exponent is of order $n^2$, the exponent is  of order $n$, so the model is again  of the form $(\ref{model})$.
%

\vspace{0.5em}
\noindent
{\bf Ising spin glasses.}
The Edwards-Anderson (spin glass) model on a graph $G=(V,E)$ with $n$ vertices is the probability distribution
\begin{equation}
\label{eq:spin_glass}
\pi(\sigma) = \frac{1}{Z(\beta)} \exp \left\{ \beta \sum_{v \sim w} J_{vw} \sigma_v \sigma_w \right\}
\end{equation}
on $\mathcal{X} := \{-1,+1\}^V$. Again, the parameter $\beta \geq 0$ is the inverse temperature, and $Z(\beta)$ is the normalizing constant. The sum is over all edges 
of the graph $G$. In contrast to the Ising model, here we have a separate interaction constant $J_{vw}$ on each edge $vw \in E$ of the graph. 
If $J_{vw} > 0$, the interaction is \emph{ferromagnetic}, meaning that the spins $\sigma_v$ and $\sigma_w$ like to align under $\pi$, whereas if $J_{vw} < 0$, the interaction is \emph{antiferromagnetic}, meaning that the spins $\sigma_v$ and $\sigma_w$ like to anti-align under $\pi$.

We take the $J_{vw}$ to be iid Rademacher random variables, taking the values $\pm 1$ with probability $1/2$, independently over the edges $vw \in E$ of the graph.
A characteristic property of these models is \emph{frustration}, meaning that not all ``constraints'' imposed by the interaction constants $J_{vw}$ can be satisfied simultaneously.
For example, take the four vertices $(0,0),(0,1),(1,0),(1,1)$ in a graph  $G \subset \mathbb{Z}^2$, and take three of the interaction constants to be positive, and one negative. Start with an arbitrary vertex and assign a spin to it arbitrarily. Going around the circle, trying to satisfy all constraints, will eventually lead to ... frustration.
%
%
%
%
%
Note that this definition makes $\pi$ a random measure, but conditional on the choice of the interaction constants $J_{vw}$, we get in (\ref{eq:spin_glass}) a fixed probability distribution $\pi$ (with quenched 
interactions).
The joint presence of quenched disorder and frustration makes 
spin glasses
very hard to analyze and to (approximately) simulate \cite{SlySun2012}.
Clearly, if the graph $G$ has bounded degree, e.g. if $G$ is the lattice $\{0,1,...,m-1\}^d$, where $d$ is fixed as $m$ goes to infinity,
then the sum in the exponent of (\ref{eq:spin_glass}) is of order $n=m^d$, so the model is again of the form (\ref{model}).
For more background on all of the above models from a statistical physics perspective, see \cite{MezardMontanari2009}.

%

\vspace{0.5em}
\noindent
{\bf Exponential random graph models.}
An exponential random graph model is a probability distribution 
\begin{equation}
\label{eq:ERGM}
\pi(G) = \frac{1}{Z(\beta)} \exp \left\{ \sum_{i=1}^k \beta_i T_i(G) \right\}
\end{equation} 
on the space of simple graphs with $\nu$ vertices. Here $\beta = (\beta_1,...,\beta_k)$ is a vector of real-valued parameters, and $(T_1,...,T_k)$ is the sufficient statistic. A concrete example from \cite{ChatterjeeDiaconis2013} is
$$
\pi(G) = \frac{1}{Z(\beta_1,\beta_2)} \exp \left\{ 2 \beta_1 E + 6 \beta_2 \frac{\Delta}{\nu} \right\} ,
$$
where $E=E(G)$ is the number of edges of the graph $G$, and $\Delta=\Delta(G)$ is the number of triangles of $G$. The scaling ensures nontrivial limits. (Without proper scaling, almost all graphs are empty or full in the large $\nu$ limit.)
Since both $E$ and $\Delta/\nu$ are of order $n := {\nu \choose 2}$, this model is also of the form (\ref{model}). 
Assuming proper scaling, the same is true for the model (\ref{eq:ERGM}) in general, as long as the number $k$ of statistics $T_i$ does not depend on $n$.
For  background on these models and pointers to the literature, see \cite{ChatterjeeDiaconis2013,GoldenbergZhengFienbergEtAl2010}.

Note that the space of all simple graphs $G$ on $\nu$ vertices is in natural one-to-one correspondence with $\mathcal{X} := \{0,1\}^n$. To see this, put all $n = {\nu \choose 2}$ potential edges in some arbitrary but fixed order. Then, for $x \in \mathcal{X}$, if the $i^{th}$ coordinate of $x$ is zero, the $i^{th}$ edge is absent in $G$. If the $i^{th}$ coordinate of $x$ is one, the $i^{th}$ edge is present in $G$. In particular, this shows that the the support of $\pi$ is simple, in the sense that it is easy to simulate the uniform distribution on it.
%
Note that in this case, the 
size parameter (dimension) is $n = \OO{\nu^2}$, where $\nu$ is the number of vertices of the graphs $G$.

\vspace{0.5em}
It is easy to see that the state spaces in all of the above examples satisfy assumption (A)(1). For a formal proof, consider, for example, the Potts model, where $\mathcal{X} := [q]^n$ for some $q \in \mathbb{N}$. Let $P^{(0)}$ be the transition kernel of the Gibbs sampler for the uniform distribution on $\mathcal{X}$.
We can use the following well known  coupling for this process: Suppose we are currently at $(V_t,W_t)$.
Draw $i \in [n]$ and $B \in [q]$ uniformly at random, independent of each other (and of all previous choices). Then set the $i^{th}$ coordinates of 
$V_{t+1}$ and $W_{t+1}$ 
both to $B$, and leave the other coordinates unchanged. Note that in this coupling, for all coordinates $i \in [n]$, we have
$$
V_s^{(i)}=W_s^{(i)} \;\;\;\; \Rightarrow \;\;\;\; V_t^{(i)}=W_t^{(i)} \; \mbox{ for all } \; t \geq s .
$$
Let $\tau_0$ be the first time that all $n$ coordinates have been chosen at least once in this coupling. By coupon collecting, if $t \geq \lceil n \log n + cn \rceil$, then
$$
P_{x,y} \left\{ V_t \neq W_t \right\} \leq P \{ \tau_0 > t \} \leq e^{-c}.
$$
See, for instance, \cite{LevinPeresWilmer2009}, Proposition 2.4 on page 23.
Here and throughout the paper, subscripts to the probability measure indicate the starting states (here: $V_0=x,W_0=y$) of the process.
So if we define
$C(\epsilon) := 3 + 2\log (1/\epsilon)$, 
then we get
$$
P_{x,y} \left\{V_t \neq W_t \right\} \leq  \frac{\epsilon}{n+1} 
$$
for all  $t \geq C(\epsilon) n \log n$, as required, since 
$C(\epsilon) n \log n \geq  \big\lceil n \left[ \log (n) + \log ((n+1)/\epsilon) \right] \big\rceil$ for all $n \geq 2$. $\Box$

\section{The algorithm}
\label{sec:Alg}
To approximately simulate from a distribution $\pi$ of the form (\ref{model}),
we use a version of the interacting tempering algorithm, see \cite[section 3 on page 3274]{FortMoulinesPriouret2011}, and also \cite{AndrieuJasraDoucetEtAl2011}.
Define tempered versions of $\pi$ by
$$
\pi_j(x) := 
\frac{1}{Z(j,\beta)} \exp\left\{ j \, \beta \, S(x) \right\}, \;\;\;\;\;\; x \in \mathcal{X},
$$
for $j=0,1,...,n$. Note that $\pi_n$ is equal to the distribution of interest $\pi$, while $\pi_0$ is the uniform distribution on $\mathcal{X}$. The distributions $\pi_j$ ``interpolate'' between $\pi_0$ and $\pi_n$,
where we use a temperature ladder with $n+1$ temperatures such that the inverse temperatures $\beta j/n$ are equally spaced over the interval from zero to $\beta$.

Let $P^{(0)}$ denote the transition kernel from assumption (A)(1) targeting the uniform distribution $\pi_0$.
For $j=1,...,n$, let $P^{(j)}$ be the transition kernel of some (any) local move Markov chain (that we can simulate) with unique stationary distribution $\pi_j$. The particular choice of the kernels $P^{(j)}$ for $j=1,...,n$ will not affect any of our results.
For concreteness, let  $P^{(j)}$ be the lazy random walk Metropolis algorithm for $\pi_j$. 
To specify this Markov chain, we first have to define an (arbitrary)  connected graph $G'$ with vertex set $\mathcal{X}$. 
The chain then evolves as follows. 
Given we are currently at state $Z_t^{(j)}$, we first flip a fair coin. If it comes up heads, we stay where we are, setting $Z_{t+1}^{(j)} = Z_t^{(j)}$. If it comes up tails, we select one of the 
neighbors $Y$ of $Z_t^{(j)}$ (in the graph $G'$) uniformly at random. 
Then, with probability 
$1 \wedge \frac{\pi_j(Y) / N(Y) }{ \pi_j(Z_t^{(j)}) / N(Z_t^{(j)})}$,
we accept the proposal and set $Z_{t+1}^{(j)} = Y$. With the remaining probability, we reject the proposal and set $Z_{t+1}^{(j)} = Z_t^{(j)}$. Here, $N(x)$ is the number of neighbors of state $x \in \mathcal{X}$ in the graph $G'$.

\vspace{0.5em}
%
%
{\bf The algorithm:} In our setting, the interacting tempering algorithm specifies a process
$(X_t) = \left(X_t^{(0)}, X_t^{(1)}, ...., X_t^{(n)}\right)$ on $\mathcal{X}^{n+1}$, started at some state $X_0 \in \mathcal{X}^{n+1}$. The component $X_t^{(j)}$ will target the distribution $\pi_j$.
The two tuning parameters of our version of the algorithm are the probability of interaction $v \in (0,1)$, and an error parameter $\epsilon > 0$.
Let $\lambda := v e^{-\beta D}$, where $\beta,D$ are the constants from assumption (A) respectively representation (\ref{model}).
Let $G_0 := G_0(\epsilon) :=  C(\epsilon) n \log (n)$, 
where the constant $C(\epsilon)$ comes from assumption (A)(1), and let
$G := G(\epsilon,\lambda) := \left\lceil \log \left( \frac{n+1}{\epsilon} \right) \Big/ \log \left( \frac{1}{1-\lambda} \right) \right\rceil$. 
Define $s_0 := 0, t_0 := G_0$ and define $s_j := G_0 + (j-1)G, t_j := G_0+jG$, for $j=1,...,n$. 

Conditional on the history $(X_0,X_1...,X_t)$ of the process so far, at the time step $t \rightarrow t+1$, the process $(X_t)$ evolves as follows. 
We let $X_{t+1}^{(0)}$ be a draw from $P^{(0)}(X_t^{(0)},\cdot)$. That is, $(X_t^{(0)})$ evolves according to our 
Markov chain $P^{(0)}$ for the uniform distribution $\pi_0$.
The components $X_t^{(j)}$, for $j=1,...,n$, evolve as follows: 
If $t < s_j$, we stay where we are and set $X_{t+1}^{(j)} := X_{t}^{(j)}$. 
If $t \geq s_j$, we flip a coin with probability of heads equal to $v$.
If it comes up tails,
we let $X_t^{(j)}$ evolve according to our local move chain and draw $X_{t+1}^{(j)}$ from $P^{(j)}(X_t^{(j)},\cdot)$. 
If it  comes up heads, we do the following.
First, we draw a proposal $Y$ from the empirical distribution of $\left(X_{t_{j-1}}^{(j-1)},...,X_{t}^{(j-1)}\right)$.
Then, we accept this proposal and set $X_{t+1}^{(j)}:=Y$ with probability $1 \wedge a_j(X_t^{(j)},Y)$, where for $x,y \in \mathcal{X}$,
\begin{equation}
\label{accept_prob_IT}
a_j(x,y) := \frac{\pi_j(y) \, \pi_{j-1}(x)}{\pi_j(x) \, \pi_{j-1}(y)}.
\end{equation}
With the remaining probability, we reject the proposal and stay where we are, setting $X_{t+1}^{(j)} := X_t^{(j)}$.
Unless otherwise mentioned, all random choices in the algorithm are understood to be made independent of all previous choices. $\Box$

\vspace{0.5em}
\begin{remark}
Note that in this algorithm, $s_j$ is the time when coordinate $X^{(j)}$ starts evolving, while 
$t_j$ 
is the time when we start collecting the history of $X^{(j)}$ to be used as proposals for cross-temperature moves in coordinate $j+1$.
That is, we allow for a burn-in of $G_0$ steps before we start collecting the history of the process $(X_t^{(0)})$ and start running $(X_t^{(1)})$. Similarly, for $j=1,...,n-1$, after we start running the process $(X_t^{(j)})$, we allow for a burn-in of $G$ steps before we start collecting its history and start running $(X_t^{(j+1)})$.
\end{remark}

\begin{remark}
The choice of the local move transition kernels $P^{(j)}$, for $j=1,...,n$, does not affect any of our results, since these results are based purely on the cross-temperature moves of the algorithm.
Only the local move kernel $P^{(0)}$ (targeting the uniform distribution on $\mathcal{X}$) from assumption (A)(1) affects the burn-in $G_0$ for coordinate zero.
\end{remark}

\begin{remark}
Note that from (\ref{accept_prob_IT}) and (\ref{model}), we get for any states $x,y \in \mathcal{X}$, 
\begin{eqnarray}
\label{accept_prob_bd}
a_j(x,y) & = & \exp\{ j \beta [S(y) - S(x)] - (j-1) \beta [S(y) - S(x)] \}  \nonumber \\
& = & \exp \{ \beta [S(y) - S(x)]  \} \\
& \geq & \exp\{ -D \beta \} . \nonumber
\end{eqnarray}
Therefore, the acceptance probabilities for cross-temperature moves from $j-1$ to $j$ are bounded away from zero (uniformly in $j$ and $n$). This was the reason for choosing $n+1$ inverse temperatures $j\beta/n$, where $j=0,1,...,n$.
Also note that the acceptance probabilities in (\ref{accept_prob_IT}) correspond to the ones we would get in the Metropolis Hastings algorithm
\emph{if the proposal $Y$ would be an independent draw from $\pi_{j-1}$.} In the actual algorithm, we \emph{approximate} this independent draw from $\pi_{j-1}$ with a draw from the empirical distribution of $\left(X_{t_{j-1}}^{(j-1)},...,X_{t}^{(j-1)}\right)$. The idea is that if the process $(X_t^{(j-1)})$ has converged (approximately) to $\pi_{j-1}$ 
by time $t_{j-1}$,
and if its mixing time is small compared to $t$,
then this should be a good approximation.
\end{remark}

\section{Main result}
\label{sec:Res}
For random elements $X,Y$, we write $\mathcal{D}(X)$ for the distribution (i.e., the law) of $X$, and we write $\mathcal{D}(X \, | \, Y)$ for the conditional distribution of $X$ given $Y$.
Subscripts indicate starting distributions (respectively starting states). For example, $\mathcal{D}_{\mu}(X_t)$ (respectively $\mathcal{D}_x(X_t)$) denotes the distribution of the interacting tempering process $(X_t)$ at time $t$, when started in distribution $\mu$ on $\mathcal{X}^{n+1}$ (respectively in state $x \in \mathcal{X}^{n+1}$).
The following statement is our main result.
\begin{theorem}
\label{thm:LOM}
For any probability distribution $\pi$ that satisfies assumption (A),
the interacting tempering process, as defined above,
forgets its starting distribution
after $G_0+nG$ steps. That is, for any $\epsilon > 0$, and for any starting 
distributions $\mu$ and $\nu$ on $\mathcal{X}^{n+1}$, the total variation distance after 
$t \geq G_0 + nG$ steps  of the algorithm (with error parameter $\epsilon$) satisfies
$$
|| \mathcal{D}_{\mu} (X_t) - \mathcal{D}_{\nu} (X_t) || \leq \epsilon. 
$$
Here, we have $G_0 := C(\epsilon) n \log n$, $v \in (0,1)$ is the probability of interaction, 
$\lambda := ve^{-\beta D}$,  and
$$
G   := \left\lceil \log \left( \frac{n+1}{\epsilon} \right) \Big/ \log \left( \frac{1}{1-\lambda} \right) \right\rceil ,
$$
where the constants  $C(\epsilon),D,\beta$ are from assumption (A) respectively representation (1).
\end{theorem}

\begin{remark}
The theorem shows that the interacting tempering algorithm for any target distribution $\pi$ that satisfies assumption (A) forgets its starting distribution in $G_0 + nG = \OO{n \log n}$ steps.
Here, an update of the process from $X_t$ to $X_{t+1}$ is counted as one step. Since one such step generally  involves updating all $n+1$ coordinates of $X_t = \left( X_t^{(0)},X_t^{(1)}, ..., X_t^{(n)} \right)$,
the computational effort to forget the starting distribution is of order 
$(n+1) \times (G_0 + nG) = \OO{n^2 \log n}$.  
\end{remark}

\begin{remark}
Since the interacting tempering process is generally not Markovian, the theorem says nothing about the (more interesting) question of how long it takes for the process to converge to its limiting distribution
$\Pi$ on $\mathcal{X}^{n+1}$.
However, the result here may be seen as a stepping stone towards such quantitative, non-asymptotic convergence rates, since, for the sake of bounding such rates, our result allows us to start the interacting tempering process in its limiting distribution 
$\Pi$.  
To see this, note that by the triangle inequality, for any starting distribution $\mu$ on $\mathcal{X}^{n+1}$,
\begin{equation}
\label{eq:ConvRates}
|| \mathcal{D}_{\mu}(X_t) - \Pi || \leq || \mathcal{D}_{\mu}(X_t) - \mathcal{D}_{\Pi}(X_t) || + || \mathcal{D}_{\Pi}(X_t) - \Pi ||.
\end{equation}
Our result gives an upper bound of $\epsilon$ for the first term on the right hand side above.
Therefore, to get rates of convergence, it remains to bound the second term on the right hand side of (\ref{eq:ConvRates}). That is, we may assume that the process starts in its limiting distribution $\Pi$. 
(This is sometimes called a \emph{warm start}.)
However, since the transition rule for $(X_t)$ does not preserve $\Pi$, bounding this remaining term
on the right hand side of (\ref{eq:ConvRates})
will generally 
not be easy.
\end{remark}

\vspace{0.5em}
\noindent {\bf Proof of Theorem \ref{thm:LOM}.} By repeated application of the triangle inequality for the total variation norm, we get
$$
|| \mathcal{D}_{\mu}(X_t) - \mathcal{D}_{\nu}(X_t) || \leq \sup_{x,y} || \mathcal{D}_x(X_t) -\mathcal{D}_y(X_t) ||,
$$
so it's enough to bound the right hand side above.
Fix any $\epsilon >0$ and any starting states $x,y \in \mathcal{X}^{n+1}$.
Recall the definition of the times $s_j$ (when coordinate $X^{(j)}$ starts evolving) and $t_j$ (when we start collecting the history of $X^{(j)}$) from the specification of the algorithm in section \ref{sec:Alg}.
We will prove the theorem by constructing a coupling $(X_t,Y_t)$ of two versions of the interacting tempering process $(X_t)$, one started at $X_0=x$, the other started at $Y_0=y$. 
By the coupling inequality, it will then be enough to show that for all $t \geq t_n$, 
\begin{equation}
\label{coupling_inequality}
P_{x,y} \{ X_t \neq Y_t \} \leq \epsilon.
\end{equation}
For any  $k \leq l$, write $X_{k:l} := (X_s)_{s=k,...,l}$ 
and $X_{k:l}^{(j)} := (X_s^{(j)})_{s=k,...,l}$ 
for the history of the entire process
(respectively, of component $j$) from steps $k$ to $l$.
Analogously, for any $u \leq v$ and $k \leq l$,  write $X_t^{(u:v)}:=(X_t^{(u)},...,X_t^{(v)})$ 
and $X_{k:l}^{(u:v)} := (X_s^{(u)},...,X_s^{(v)})_{s=k...,l}$, 
for coordinates $u$ to $v$ at time $t$ (respectively, from time $k$ to $l$).
Note that, by construction of the algorithm, the one step transition probabilities for higher temperature coordinates do not depend on the history of the lower temperature coordinates, once we condition on the history of those higher temperature coordinates. That is, for any $j=0,1,...,n$, we get
\begin{equation}
\label{MarkovProperty}
\mathcal{D} \left( X_{t+1}^{(0:j)} \, | \, X_{0:t}^{(0:n)} \right) = \mathcal{D} \left( X_{t+1}^{(0:j)} \, | \, X_{0:t}^{(0:j)} \right).
\end{equation}
This allows us to 
work 
by induction on $j$.
For $j=0$, we note that, marginally, the process $(X_t^{(0)})$ is a time homogeneous Markov chain with transition kernel $P^{(0)}$.
To define a coupling $(X_t^{(0)},Y_t^{(0)})$ of two versions of this Markov chain, 
we simply use the coupling that exists by assumption (A)(1). (See the end of section \ref{sec:examples} for an explicit construction of such a coupling 
for the state space $\mathcal{X} = [q]^n$
of the Potts model.)
Next we define the coupling $(X_t,Y_t)$ for  coordinates $j = 1,...,n$. 
It will be enough to specify how to do one step $t \rightarrow t+1$.
Let $A_j := \{ X_{t_j}^{(j)} = Y_{t_j}^{(j)} \}$, and let $B_j := \bigcap_{i=0}^j A_i$, for $j=0,1,...,n$. 
Suppose our coupling is already specified for coordinates $i=0,1,...,j-1$, and suppose 
the history so far is
$(X_{0:t}^{(0:j)},Y_{0:t}^{(0:j)})$. From the induction hypothesis, 
we can draw $(X_{t+1}^{(0:(j-1))},Y_{t+1}^{(0:(j-1))})$ according to our coupling as already defined. It remains to specify how we draw $(X_{t+1}^{(j)},Y_{t+1}^{(j)})$, conditional on the history
$(X_{0:t}^{(0:j)},Y_{0:t}^{(0:j)})$.
(By construction of the algorithm, the transition from time $t$ to time $t+1$ only depends on the history of the process up to time $t$,
and this will also be true for our coupling.)
If $t < s_j$, we don't move and set $(X_{t+1}^{(j)},Y_{t+1}^{(j)}) := (X_t^{(j)},Y_t^{(j)})$.
If $t \geq s_j$, we proceed as follows.
On the complement of the event $B_{j-1}$, we let both processes $(X_t^{(j)})$ and $(Y_t^{(j)})$ evolve independently according to their respective transition rules. On the event $B_{j-1}$, 
we do the following.
Flip a coin with probability of heads equal to $v$. 
\begin{itemize}
\item If it comes up heads, we attempt a cross-temperature move, and do the following: Let $Z'$ be a draw from the empirical distribution of the history $(X_{t_{j-1}}^{(j-1)},...,X_{t}^{(j-1)})$, and let 
$U'$ be an (independent) $Uniform(0,1)$ random variable. Then set
\begin{eqnarray*}
X_{t+1}^{(j)} & := & \begin{cases} Z' & : \mbox{ if } U' < a_j(X_t^{(j)},Z') , \\
			 				   X_{t}^{(j)} & : \mbox{ otherwise, }
				 \end{cases} \\
Y_{t+1}^{(j)} & := & \begin{cases} Z' & : \mbox{ if } U' < a_j(Y_t^{(j)},Z') , \\
							   Y_{t}^{(j)} & : \mbox{ otherwise.  }
				 \end{cases} 
\end{eqnarray*}
\item If it comes up tails, we make local moves, and do the following:
Draw $X_{t+1}^{(j)} \sim P^{(j)}(X_t^{(j)}, \cdot)$ and (independently) draw $Z'' \sim P^{(j)}(Y_t^{(j)}, \cdot)$.
Then set
$$
Y_{t+1}^{(j)} := \begin{cases} X_{t+1}^{(j)} & : \mbox{ if } X_t^{(j)}=Y_t^{(j)} , \\
							   Z'' & : \mbox{ otherwise.  }
				 \end{cases} 
$$
\end{itemize}
{\bf Claim: } 
In this coupling, for all $j=1,...,n$ and $s \geq 0$, on the event $B_{j-1}$, we have 
\begin{equation}
\label{XeqYforever}
X_{s}^{(j)} = Y_{s}^{(j)}  \;\;\;\;\; \Rightarrow \;\;\;\;\; X_t^{(j)} = Y_t^{(j)} \; \mbox{ for all } t \geq s .
\end{equation}
Furthermore, (\ref{XeqYforever}) also holds for $j=0$ and all $s \geq 0$.
\\
{\bf Proof of Claim:} This follows by induction on $j$ from the construction of the coupling:
We may assume that (\ref{XeqYforever}) is true for $j=0$, since the coupling from  assumption (A)(1) is Markovian. (So if we have $X_s^{(0)}=Y_s^{(0)}$ for some $s\in \mathbb{N}$, we can always change the coupling to ensure that 
$X_t^{(0)}=Y_t^{(0)}$ for all $t \geq s$.) 
The induction step from $j-1$ to $j$ is an immediate consequence of the construction of the coupling. $\Box$

To see that the algorithm specified above is a valid coupling of the two copies of our process, note that on the event $B_{j-1}$, by the Claim we have agreement of the histories $(X_{t_{j-1}}^{(j-1)},...,X_{t}^{(j-1)})$ and $(Y_{t_{j-1}}^{(j-1)},...,Y_{t}^{(j-1)})$.
We may therefore make \emph{one common draw} from this \emph{joint history} when proposing 
cross-temperature moves, as we did in the construction above. 

Note that on the event $B_n$, we get by the Claim that $X_t=Y_t$ for all $t \geq t_n$.
Therefore, to establish (\ref{coupling_inequality}) for our coupling $(X_t,Y_t)$, it will be enough to show the inequality
\begin{equation}
\label{product_eps}
P_{x,y} (B_n) = P_{x,y} \Big( \bigcap_{j=0}^n A_j \Big) = P_{x,y} (A_0) \, \prod_{j=1}^n P_{x,y} \Big( A_j \, | \, B_{j-1} \Big) \geq 1-\epsilon.
\end{equation}
To establish (\ref{product_eps}), it suffices to show that
\begin{equation}
\label{product_eps2}
P_{x,y} (A_0) \geq 1-\frac{\epsilon}{n+1} \;\; \mbox{ and }  \;\;\; P_{x,y} \Big( A_j \, | \, B_{j-1} \Big) \geq 1-\frac{\epsilon}{n+1} \; \mbox{ for all } j=1,...,n,
\end{equation}
since this implies
$$
P_{x,y} \Big( \bigcap_{j=0}^n A_j \Big) \geq \Big( 1-\frac{\epsilon}{n+1} \Big)^{\!n+1} \geq 1-\epsilon.
$$
To see the last inequality above, we can use calculus to show that the function that maps $\epsilon$ to $(1-\epsilon/(n+1))^{n+1} - (1-\epsilon)$ is nonnegative. 
To see that our coupling satisfies 
(\ref{product_eps2}),
consider the $j^{th}$ coordinate of the coupling. 
For $j=0$, this follows by assumption (A)(1), so suppose $j \geq 1$.
On the event $B_{j-1}$, 
for cross-temperature moves in our coupling we always use the \emph{same} proposal $Z'$ for both $X_t^{(j)}$ and $Y_t^{(j)}$.
Therefore, we get $X_t^{(j)}=Y_t^{(j)}$ as soon as we accept a cross-temperature move \emph{in both coordinates at the same time}. Starting at $s_{j}$, 
the time when the $j^{th}$ component starts evolving,
we see from (\ref{accept_prob_bd}) that the time until this happens is stochastically dominated by a Geometric random variable with success probability $\lambda := ve^{-D\beta}$. This means we get
$$
P_{x,y} \Big( A_j^c \, \Big| \, B_{j-1} \Big) \leq (1-\lambda)^G \leq \frac{\epsilon}{n+1} ,
$$
where the last inequality above follows by the definition of
$$
G := 
\left\lceil \log \left( \frac{n+1}{\epsilon} \right) \Big/ \log \left( \frac{1}{1-\lambda} \right) \right\rceil .
$$
This establishes (\ref{product_eps2}), finishing the proof of the theorem. $\Box$

\section{A cautionary example}
\label{sec:NIH}

In this section we exhibit an example of a probability distribution $\pi$ on $\mathcal{X}:=\{0,1\}^n$ with the following two properties: 
%
\begin{enumerate}
\item The distribution $\pi$ satisfies our assumption (A) of simple support and exponentially bounded likelihood ratios.
Consequently, by Theorem \ref{thm:LOM}, the interacting tempering algorithm for $\pi$ forgets its starting distribution in order $n \log n$ steps. 
\item It takes at least an exponential (in $n$) number of steps 
for the $n^{\textrm{th}}$ coordinate of the interacting tempering algorithm to get close to its limiting distribution $\pi$. This property is often called torpid mixing.
\end{enumerate}
More specifically, we have the following result.
\begin{theorem}
\label{thm:NIH}
For every $\epsilon>0$, there exists a 
probability distribution $\pi_{\epsilon}$ 
on $\mathcal{X}:=\{0,1\}^n$
that satisfies assumption (A), and such that the following is true.
For every starting distribution $\mu$,
the interacting tempering process $(X_t)$ for $\pi_{\epsilon}$
(as defined in section \ref{sec:Alg}, with error parameter $\epsilon/4$),
satisfies
$$
|| \mathcal{D}_{\mu}\left(X_t^{(n)}\right) - \pi_{\epsilon} || \geq 1-\epsilon 
$$
for all $t \in [t_n',\epsilon (n+1)^{-1} 2^{n-2} - 1]$, where $t_n' := G_0(\epsilon/4) + n G(\epsilon/4)$,
and
$$
G_0(\epsilon/4) := \Big\lceil n \left[ \log (n) + \log \left(4(n+1)/\epsilon\right) \right] \Big\rceil, \;\;\;
G(\epsilon/4)   := \left\lceil \log \left( \frac{n+1}{\epsilon/4} \right) \Big/ \log \left( \frac{1}{1-\lambda} \right) \right\rceil , \;\;\;
\lambda := v \epsilon / 4.
$$ 
\end{theorem}

\begin{remark}
The theorem shows that for any starting distribution $\mu$, it takes at least an exponential (in $n$) number of steps for this interacting tempering algorithm to get close to its limiting distribution $\pi_{\epsilon}$.
Since at time $t_n'$, the coordinate of interest $\left(X_t^{(n)}\right)$ has only made $G(\epsilon/4)=\OO{\log n}$ moves, it is not a real restriction that the theorem remains silent about times $t<t_n'$.
\end{remark}
\vspace{0.5em}
We will use the following family of distributions to establish Theorem \ref{thm:NIH}. Fix $n \in \mathbb{N}$ and chose an arbitrary state $z \in \mathcal{X} := \{0,1\}^n$.
For each $\delta \in (0,1)$, define
\begin{equation}
\label{eq:NeedleDelta}
\tilde{\pi}(x) := 
\begin{cases} 2^n/\delta & : \; x=z,\\ 1 & : \; x \in \mathcal{X}, x \neq z . \end{cases}
\end{equation}
Then let $\pi(x) := 
\tilde{\pi}(x)/Z$, where 
$Z = Z(\delta) := \sum_x \tilde{\pi}(x) = 2^n(1+\delta^{-1})-1$ is the normalizing constant.
An easy calculation shows  that $\pi(z) \geq 1-\delta$.
For the distribution $\pi$, the state $z$ can be thought of as a ``needle'' in the (exponential size) ``haystack'' $\mathcal{X}$.

In the interacting tempering algorithm, as defined in section \ref{sec:Alg}, we use tempered versions $\pi_j(x) \propto \pi(x)^{j/n}$ for $j=0,1,...n$, so that $\pi_0$ is the uniform distribution $U$ on $\mathcal{X}$, and $\pi_n$ is the target distribution $\pi$.
As local move transition kernels $P^{(j)}$, we use 
the lazy random walk Metropolis algorithm as specified in the definition of the algorithm in section \ref{sec:Alg}. For this, we define the required graph $G'$ on $\mathcal{X}=\{0,1\}^n$ by saying that $x,y \in \mathcal{X}$ are neighbors in $G'$ iff they differ in exactly one coordinate.

Note that for all $x,y \in \mathcal{X}$, we have $\pi(x)/\pi(y) = \tilde{\pi}(x) / \tilde{\pi}(y) \leq 2^n/\delta$, so $\pi$ satisfies assumption (A) with $D := \log(2/\delta)$, and $\beta=1$ in representation (\ref{model}).
Consequently, by Theorem \ref{thm:LOM}, the interacting tempering algorithm for $\pi$ forgets its starting distribution in order $n \log n$ steps.
To prove  Theorem \ref{thm:NIH}, we have to show that for any starting distribution $\mu$, it takes at least an exponential number of steps for the interacting tempering process to get close to $\pi$. To do this, we first analyze the case where the starting state for each coordinate is chosen from the uniform distribution on $\mathcal{X}$.
The general case will then be reduced to this special case.
\begin{prop}
\label{prop:NeedleUnif}
Fix any $\epsilon,\delta>0$, and consider the interacting tempering process $(X_t)$ targeting 
the distribution 
$\pi_{\delta}$
as defined in (\ref{eq:NeedleDelta}), 
with arbitrary choice of interaction probability and error parameter.
Then, for all $t \leq \epsilon (n+1)^{-1} 2^{n-2} -1$, we get
$$
|| \mathcal{D}_{\nu}\left(X_t^{(n)}\right) - \pi_{\delta} || \geq 1 - \delta - \epsilon/4 , 
$$
as long as all $n+1$ marginals of the starting distribution $\nu$ are uniform on $\mathcal{X}$.
\end{prop}
The proof of this proposition relies on the following result.
\begin{lemma}
\label{lemma:NeedleUnif}
Let $P_{\nu}^{IT(U)}$ denote the probability measure governing the interacting tempering process for the uniform target distribution $U$, when started in a distribution $\nu$ on $\mathcal{X}^{n+1}$. Then, 
for any choice of interaction probability and error parameter,
for all times $t \in \mathbb{N}$,  all coordinates $j=0,1,...,n$, and all states $x \in \mathcal{X}$, we get
$$
P_{\nu}^{IT(U)} \{ X_t^{(j)} = x \} = U(x),
$$
provided that all $n+1$ marginals of $\nu$ are equal to $U$.
That is, if we start all coordinates in the uniform distribution on $\mathcal{X}$, then marginally they will all remain in the uniform distribution  forever. 
\end{lemma}

\noindent
{\bf Proof of Lemma \ref{lemma:NeedleUnif}.} 
Since the target distribution is already uniform, tempering has no effect, and we get $\pi_j = U$ for all $j=0,1,...,n$. 
Consequently, cross-temperature moves are always accepted in the interacting tempering process targeting $U$. %
The result now follows essentially from the fact that $U$ is stationary for our local move kernels $P^{(j)}$ (lazy simple random walk),
together with the fact that the cross-temperature move times are independent of the states. 

To prove this formally, we use induction on $j$ and 
on $t$. 
Fix $t \in \mathbb{N}$. By assumption, we have $\mathcal{D}(X_0^{(j)})=U$ for all $j=0,1,...,n$.
For coordinate $j=0$, there are no cross-temperature moves, and the uniform distribution $U$ is stationary for the local move transition kernel $P^{(0)}$. This shows that $\mathcal{D}(X_s^{(0)})=U$ for all $s \leq t$.
Then suppose we already know that $\mathcal{D}(X_s^{(j-1)})=U$ for all $s \leq t$ and suppose $\mathcal{D}(X_{t-1}^{(j)})=U$. 
To show that $\mathcal{D}(X_t^{(j)})=U$, we condition on the $Ber(v)$ random variable that decides whether we make a local move or a cross-temperature move in coordinate $j$ at time $t-1$. Then we condition on the time index $i$ of the draw from the history of $(X_s^{(j-1)})$, and lastly on the state $X_{t-1}^{(j)}$ respectively $X_i^{(j-1)}$. This shows that for all $x \in \mathcal{X}$,
\begin{eqnarray*}
P(X_t^{(j)}=x) & = & v P(X_t^{(j)}=x \, | \, \mbox{cross-temp. move}) + (1-v) P(X_t^{(j)}=x \, | \, \mbox{local move}) \\
& = & v \frac{1}{t-t_{j-1}} \sum_{i=t_{j-1}}^{t-1} P(X_i^{(j-1)}=x) + (1-v) \sum_{y \in \mathcal{X}} 
P^{(j)}(y,x) P(X_{t-1}^{(j)}=y) \\
& = & v U(x) + (1-v) U(x) \\
& = & U(x).
\end{eqnarray*}
Here the last but one equality follows from the induction hypothesis and the fact that $U$ is stationary for the transition kernel $P^{(j)}$ of the local move chain for $\pi_j=U$.
For the second equality, recall that we always accept cross-temperature moves in this interacting tempering algorithm for the uniform distribution. This finishes the induction step and we are done. $\Box$

\vspace{0.5em}
\noindent
{\bf Proof of Proposition \ref{prop:NeedleUnif}.} 
Fix any $\epsilon,\delta>0$.
Also, fix any starting distribution $\nu$ on $\mathcal{X}^{n+1}$ such that all $n+1$ marginals of $\nu$ are equal to the uniform distribution $U$ on $\mathcal{X}$.
By construction of the interacting tempering algorithm for $\pi_{\delta}$, local moves are always accepted, except (possibly) if the current state is the special state $z$.
Furthermore, since we have
$S(x) := n^{-1} \log(\pi(x)/\pi_{min}) 
= \mathds{1}_{\{x=z\}} \log(2 \delta^{-1/n}),$
and
$a_j(x,y) = \exp\{ \beta [S(y)-S(x)] \}$
for all $x,y \in \mathcal{X}$ and all $j=1,...,n$,
we see that cross-temperature moves are also always accepted, except (possibly) if the current state is $z$.
Then for $j=0,1,...,n$, let
$$
\tau_j := \inf \{ t \geq 0 : X_t^{(j)}=z \}
$$
be the hitting time of state $z$ for coordinate $j$, and let $\tau := \min \{ \tau_j : j=0,1,...,n \}$ be the (overall) hitting time of state $z$ in our interacting tempering process $(X_t)$. 

Fix a time $T \in \mathbb{N}$. The key  observation is that until time $\tau$, the coordinates of our process $(X_t)$ perform just lazy simple random walk on $\mathcal{X}$, except that 
at the event times of an independent Bernoulli($v$) process, 
the next state is drawn from the empirical distribution of the history of the process one step higher up in the temperature ladder. But this corresponds exactly to the interacting tempering algorithm where the target distribution is $U$ instead of $\pi_{\delta}$. Until time $\tau$, these two interacting tempering processes follow the exact same law. Consequently, by summing over paths,
we get
\begin{eqnarray*}
P_{\nu}^{IT(\pi_{\delta})} \left\{ \tau > T \right\} & = & P_{\nu}^{IT(U)} \left\{ \tau > T \right\} \\
& = & 1 - P_{\nu}^{IT(U)} \left( \bigcup_{j=0}^n \bigcup_{t=0}^T \left\{ X_t^{(j)} = z \right\} \right) \\
& \geq & 1 - \sum_{j=0}^n \sum_{t=0}^T P_{\nu}^{IT(U)} \{ X_t^{(j)}=z \} \\
& = & 1 - \sum_{j=0}^n \sum_{t=0}^T U(z) \\
& = & 1 - (n+1) (T+1) \, 2^{-n}.
\end{eqnarray*}
The the last but one equality 
in the above 
comes from Lemma \ref{lemma:NeedleUnif}.
For $A := \mathcal{X}\backslash\{z\}$, we then get
\begin{eqnarray*}
|| \mathcal{D}_{\nu}(X_T^{(n)}) - \pi_{\delta} || & \geq & P_{\nu}\{X_T^{(n)} \in A\} - \pi_{\delta}(A) \\
& \geq & 1 - P_{\nu}\{X_T^{(n)} = z\} - \delta \\
& \geq & P^{IT(\pi_{\delta})}_{\nu}\{\tau > T\} - \delta .
\end{eqnarray*}
Here, the second inequality uses $\pi_{\delta}(z) \geq 1-\delta$, while the last inequality follows since $X_T^{(n)}=z$ implies $\tau \leq T$.
%
%
Combining this with the inequality above, we get
\begin{eqnarray*}
|| \mathcal{D}_{\nu}(X_T^{(n)}) - \pi_{\delta} || & \geq & 1 - (n+1)(T+1)2^{-n} - \delta \\
& \geq & 1 - \epsilon/4 - \delta,
\end{eqnarray*}
whenever $T \leq \epsilon (n+1)^{-1} 2^{n-2}-1$.
This finishes the proof. $\Box$


\vspace{0.5em}
\noindent
{\bf Proof of Theorem \ref{thm:NIH}.} We reduce the general case of Theorem \ref{thm:NIH} to the special case of marginally uniform distributions in  Proposition \ref{prop:NeedleUnif}.
Fix any $\epsilon>0$ and set $\delta := \epsilon/2$. 
Let $(X_t)$ be the interacting tempering algorithm with error parameter $\epsilon/4$, targeting the distribution $\pi_{\delta}$ as defined in (\ref{eq:NeedleDelta}).
Fix any starting distribution $\mu$ on $\mathcal{X}^{n+1}$, and let $\nu := \bigotimes_{j=0}^n U$ be the $(n+1)-$fold product of uniform distributions $U$ on $\mathcal{X}$.
Note that the definition of $t_n' := G_0(\epsilon/4) + n G(\epsilon/4)$ in Theorem \ref{thm:NIH} corresponds to the definition of $t_n = G_0(\epsilon) + n G(\epsilon)$ in Theorem \ref{thm:LOM}, except that we changed the error parameter of the algorithm from $\epsilon$ to $\epsilon/4$.
Consequently, by the triangle inequality for the total variation norm, 
we see that for all $t \geq t_n'$,
\begin{eqnarray*}
|| \mathcal{D}_{\nu}(X_t^{(n)}) - \pi_{\delta} || & \leq & 
|| \mathcal{D}_{\nu}(X_t^{(n)}) - \mathcal{D}_{\mu}(X_t^{(n)})  ||  + || \mathcal{D}_{\mu}(X_t^{(n)}) - \pi_{\delta} || \\
    & \leq & \epsilon/4 + || \mathcal{D}_{\mu}(X_t^{(n)}) - \pi_{\delta} || .
\end{eqnarray*}
The last inequality above comes from Theorem \ref{thm:LOM}, expressing the fact that by time $t_n'$ we will have forgotten the starting distribution in this interacting tempering process.
(By adapting the coupling at the end of section \ref{sec:examples} to the situation here, we can see that the constant $G_0(\epsilon/4)$ defined in Theorem \ref{thm:NIH} 
satisfies the assumptions of Theorem \ref{thm:LOM}.)
The special case in Proposition \ref{prop:NeedleUnif} gives  a lower bound on the left hand side above: For all $t \leq \epsilon (n+1)^{-1}2^{n-2}-1$, 
we get 
$$
|| \mathcal{D}_{\nu}(X_t^{(n)}) - \pi_{\delta} || \geq 1 - \epsilon/4 - \delta.
$$
Combining the two inequalities above finishes the proof. $\Box$

\section{Convergence diagnostics}
\label{sec:Diag}
Suppose our goal is to estimate the expectation $\mu := \sum_{x \in \mathcal{X}} h(x) \pi(x)$ of a function $h:\mathcal{X} \rightarrow \mathbb{R}$ with respect to the distribution $\pi$. Let $(Y_t)$ be a process that converges to $\pi$. For example, $(Y_t)$ could be a time-homogeneous Markov chain that is ergodic to $\pi$; or, $(Y_t)$ could be the $n^{th}$ coordinate $(X_t^{(n)})$ of our interacting tempering algorithm $(X_t)$ targeting $\pi$. Discarding a burn-in of $B$ steps to reduce bias, we could use the Monte Carlo estimator
\begin{equation}
\label{eq:Est}
\hat{\mu}_{t} 
:= \frac{1}{t} \sum_{s=B+1}^{B+t} h(Y_s) .
\end{equation}
In practical applications of MCMC methods, we rarely have rigorous bounds on the convergence rates of the  process $(Y_t)$. (But see, e.g., \cite{JonesHobert2001,JonesHobert2004,LatuszynskiMiasojedowNiemiro2013,Rosenthal1993, Rosenthal1995} for some of the rare exceptions.) Consequently, there are typically no guarantees that the resulting estimate $\hat{\mu}_t$ will be within a certain margin of error of the true value $\mu$ with high probability. 
In the absence of such  rigorous guarantees, we usually rely on \emph{convergence diagnostics}. This amounts to testing certain necessary (but not sufficient) conditions for convergence of the process $(Y_t)$ to its limiting distribution $\pi$. Consequently, only negative answers come with a guarantee: If the diagnostic tells us that the process has not converged yet, this answer will generally be reliable. On the other hand, if the diagnostic tells us that the process has indeed converged, we can never be sure whether that is true.
(The outcome could always be a false positive due to metastability.)
However, particularly if we use several diagnostics on the same process, and they all come back positive, saying that the process has converged, we can argue that this constitutes evidence (in a Popperian sense) for the hypothesis that the process has indeed converged. After all, we tried to disprove this hypothesis in several different ways, but failed to do so.

Many different convergence diagnostics have been proposed in the literature, and are in use in practical applications. For an overview, see \cite[chapter 12]{CowlesCarlin1996,GelmanShirley2011, RobertCasella2004}.
Informal diagnostic procedures often involve \emph{trace plots} of $h(Y_t)$ 
against $t$, and judging whether the resulting plots ``look stationary''. For example, a clear upward (or downward) drift over such an entire plot would be evidence against stationarity, suggesting that the process $(Y_t)$ has not yet been run long enough. Another popular informal diagnostic involves running $m$ independent copies of the process $(Y_t)$, starting from different starting points,
and then producing an overlay of  $m$ trace plots. 
If the plots for different starting points do not ``overlap'' sufficiently, the influence of the starting values of the different processes might still be present,  indicating lack of convergence. 
This idea has been formalized in a widely used diagnostic by Gelman and Rubin \cite{GelmanRubin1992}. See also \cite[section 2.1]{GelmanShirley2011,CowlesCarlin1996}.  
Their diagnostic is based on computing the variance of the simulated values $h(Y_t)$ in each of the $m$ independent processes (after discarding a burn-in), 
then averaging these $m$ within process variances, and then comparing this average to the variance of the simulated values from all the $m$ processes mixed together.
At convergence, the ratio of this ``mixture variance'' to the average within process variance 
should be (close to) one, so a value of the ratio substantially larger than one is taken as evidence that the processes have not yet converged.

Our results suggest that ``mixing between processes'', in the sense of diminishing influence of the starting values of the different processes, might often happen a lot faster in the interacting tempering algorithm, than ``mixing within processes'', in the sense of convergence to the limiting distribution, a setting reported to be unusual in non-adaptive MCMC settings \cite[page 165]{GelmanShirley2011}.
We argue that this should be taken into account when performing convergence diagnostics for adaptive MCMC algorithms like interacting tempering.

To explain this difference, recall that for a time-homogeneous Markov chain $(Y_t)$, forgetting the starting distribution is equivalent to convergence to the limiting distribution $\pi$, in the following sense.
Define
$d(t) := \sup_{x \in \mathcal{X}} || \mathcal{D}_x(Y_t) - \pi ||$ and
$\bar{d}(t) := \sup_{x,y \in \mathcal{X}} || \mathcal{D}_x(Y_t) - \mathcal{D}_y(Y_t) ||$. 
Then, it is well known that for all $t \in \mathbb{N}$,
$$
d(t) \leq \bar{d}(t) \leq 2 d(t).
$$
Here, the second inequality comes from the triangle inequality for the total variation norm, whereas the first inequality relies on the fact that $\pi$ is stationary for the transition kernel of the chain. See, for instance,  \cite[Section 4.4]{LevinPeresWilmer2009}.
Of course, it is also well known that for time inhomogeneous Markov chains (and even more so for non Markovian processes like interacting tempering), forgetting the starting distribution is  necessary, but generally no longer sufficient for convergence to the limiting distribution. 
(The first inequality above is generally no longer true.)
See, for instance, \cite[Chapter 7]{Iosifescu1980} or \cite[chapter 5]{IsaacsonMadsen1976}.
Therefore, it seems a priori clear that demonstrating forgetting of the starting distribution is a much less stringent test for convergence for an adaptive MCMC algorithm, as compared to a non-adaptive MCMC algorithm based on a time homogeneous Markov chain. 
Consequently, we would argue that diagnostics based on forgetting the starting distribution are a priori less useful for adaptive MCMC algorithms, as compared to 
non-adaptive MCMC algorithms.

We believe that our results illustrate that this difference is not just purely theoretical, but that it might be of practical relevance:
There are many distributions $\pi$ that satify our assumptions (A), but that are widely believed to be very hard to (approximately) simulate. 
For example, an Ising spin glass on a regular graph
and in the anti-ferromagnetic case (where all interaction constants are minus one)
cannot be approximately simulated in time polynomial in $n$, unless NP=RP
(a complexity theoretic assumption widely believed to be false)
 \cite{SlySun2012}.
It also seems clear that our ``needle in a haystack'' model from section \ref{sec:NIH} is exponentially hard to simulate for any algorithm that only has oracle access to the unnormalized measure $\tilde{\pi}$ (since it would take an exponential number of calls to the oracle to have a positive chance to find the ``needle'' in the large $n$ limit). However, Theorem \ref{thm:LOM} shows that with appropriate choice of temperatures and burn-in for each coordinate, the interacting tempering algorithm rapidly forgets its starting distribution in all of these models. Consequently, any convergence diagnostic that is based on comparing statistics of independent copies of the process (with different starting points) would easily be fooled by this algorithm for these models. 

In summary, we belief our results suggest caution in the use of ``between processes'' diagnostics for adaptive MCMC algorithms like interacting tempering, since
these diagnostics are usually based on demonstrating that the process has forgotten its starting distribution, which is not as reliable an indicator for convergence in adaptive MCMC algorithms as it is in the non-adaptive setting.
On the other hand, ``within process'' diagnostics (like informally checking stationarity in a trace plot, or any number of other diagnostics, see \cite{CowlesCarlin1996,GelmanShirley2011, RobertCasella2004}) should remain just as valid for adaptive MCMC algorithms like interacting tempering, as they are for non-adaptive MCMC algorithms based on time-homogeneous Markov chains.

\bibliographystyle{abbrv}
\bibliography{/Users/wbarta/GW/JabRefDB/Papers}{}

\end{document}